\documentclass%[review]
{elsarticle}

\usepackage[margin=2.5cm]{geometry}%
\usepackage{mathtools}
\usepackage{epsfig}
\usepackage{amsmath}
\usepackage{amsfonts}
\usepackage{amsthm}
\usepackage{amsbsy}
\usepackage{appendix}
\usepackage{epstopdf}
\usepackage[table]{xcolor}
\usepackage{subfloat}
\usepackage{float}
\usepackage[thinlines]{easytable}
\usepackage{array}
\usepackage{color}
\usepackage{inputenc}
\usepackage{hyperref}
\usepackage[ruled,vlined]{algorithm2e}
\usepackage{verbatim}
\usepackage{caption}
\usepackage{titlesec}

\titleformat{\paragraph}
{\normalfont\normalsize\itshape}{\theparagraph}{1em}{}
\titlespacing*{\paragraph}
{0pt}{3.25ex plus 1ex minus .2ex}{1.5ex plus .2ex}

\numberwithin{equation}{section}
\allowdisplaybreaks

\theoremstyle{plain}

\newtheorem{assump}{\bf Assumption}[section]

\newtheorem{definition}{\bf Definition}[section]

\theoremstyle{remark}
\newtheorem{remark}{\bf Remark}[section]
\theoremstyle{remark}

\journal{Journal of Computational Physics}

\begin{document}

\begin{frontmatter}

\title{Learning Nonautonomous Systems via Dynamic Mode Decomposition}
\author[myaddress1,myaddress2]{Hannah Lu\fnref{email1}}
\fntext[email1]{email: \texttt{hannahlu{@}mit.edu}}
\author[mainaddress]{Daniel M. Tartakovsky\corref{mycorrespondingauthor}}
\cortext[mycorrespondingauthor]{Corresponding author}
\ead{tartakovsky@stanford.edu}
\address[myaddress1]{Department of Aeronautics and Astronautics, Massachusetts Institute of Technology, Cambridge, MA 02139, USA}
\address[myaddress2]{Department of Civil and Environmental Engineering, Massachusetts Institute of Technology, Cambridge, MA 02139, USA}
\address[mainaddress]{Department of Energy Science and Engineering, Stanford University, Stanford, CA 94305, USA}

\begin{abstract}
We present a data-driven learning approach for unknown nonautonomous dynamical systems with time-dependent inputs based on dynamic mode decomposition (DMD). To circumvent the difficulty of approximating the time-dependent Koopman operators for nonautonomous systems, a modified system derived from local parameterization of the external time-dependent inputs is employed as an approximation to the original nonautonomous system. The modified system comprises a sequence of local parametric systems, which can be well approximated by a parametric surrogate model using our previously proposed framework for dimension reduction and interpolation in parameter space (DRIPS). The offline step of DRIPS relies on DMD to build a linear surrogate model, endowed with reduced-order bases (ROBs), for the observables mapped from training data. Then the offline step constructs a sequence of iterative parametric surrogate models from interpolations on suitable manifolds, where the target/test parameter points are specified by the local parameterization of the test external time-dependent inputs. We present a number of numerical examples to demonstrate the robustness of our method and compare its performance with deep neural networks in the same settings.

\end{abstract}


\end{frontmatter}

%\linenumbers

\section{Introduction}
%% growing interest in data-driven learning
With the rapid advances in modern machine learning algorithms and the increasing availability of observational data, there have been growing developments in data-driven learning of complex systems. Various regression techniques~\cite{schmidt2009distilling,schaeffer2017learning} use a proposed dictionary, comprising plausible spatial and/or temporal derivatives of a state variable, to ``discover" the governing equations from observations. Examples of such numerical methods include sparse identification of nonlinear dynamical systems~\cite{brunton2016discovering}, Gaussian process regression~\cite{raissi2017machine}, DMD with dictionary learning~\cite{li2017extended,korda2018convergence,williams2015data} and deep neural networks (DNNs)~\cite{rico1993continuous,raissi2018deep,rudy2019deep,bakarji2021data}. Alternatively, a surrogate (aka reduced-order) model can be constructed from data in a dictionary/equation-free fashion. In this context, DMD can be used to construct an optimal linear approximation model for the unknown system~\cite{schmid2010dynamic} and to learn the unknown dynamics of chosen observables, rather than of the state itself~\cite{tu2013dynamic}. The latter is accomplished by utilizing the Koopman operator theory~\cite{koopman1931hamiltonian} in order to construct linear models on the observable space, instead of seeking for nonlinear models on the state space~\cite{brunton2017chaos}. Likewise, DNN can be used to build nonlinear surrogate models for ODEs~\cite{qin2019data,rudy2019deep} and PDEs~\cite{raissi2019physics,long2019pde,wu2020data}.

%% difficulties in nonautonomous system and literature
Nonautonomous systems have increased complexity as they describe the dynamic behaviors of physical laws due to temporal variations of system parameters, sources term, boundary conditions and etc. Learning such systems  from data is extremely challenging as the model has to differentiate the forces driven by internal dynamics from the forces driven by the external factors or inputs that vary with time. Several studies have explored the application of data-driven methods for learning nonautonomous systems in the context of control and optimization~\cite{proctor2016dynamic,proctor2018generalizing,proctor2016including,brunton2016sparse}. These methods, however, rely on a fixed structure of time-independent Koopman operator and separated time-dependent forcing/control term. More recently, a time-dependent Koopman operator framework was introduced to nonautonomous dynamical systems in~\cite{mezic2016koopman}. This extension introduces time-dependent eigenfuctions, eigenvalues and modes of the nonautonomous Koopman operator. Further developments along this line of research include multi-resolution DMD with multiple time scale decomposition~\cite{kutz2016multiresolution}, spatiotemporal pattern extraction~\cite{giannakis2019data}, delay-coordinate maps~\cite{das2019delay}. However, these methods only apply to nonautonomous systems within certain categories (e.g., periodic and quasi-periodic) and/or require special analytic tools (e.g., rescaling) and particular knowledge about the system properties. Meanwhile, there remains numerical challenges in accurate and efficient computations of the nonautonomous Koopman operator spectrum. Online DMD and weighted DMD~\cite{zhang2019online} are developed to recover the approximations of time-dependent Koopman eigenvalues and eigenfunctions. These methods require large amount  of data to capture the variations in the Koopman operator in each time window and an intrinsic error is introduced when the DMD method with the state observables is applied, which is further studied and analyzed in~\cite{macesic2018koopman}.

%% our method
Instead of approximating the time-dependent Koopman operator, our framework transforms learning of nonautonomous systems into learning of locally parameterized systems. This transformation was first used in~\cite{qin2021data} for recovering unknown nonautonomous dynamical systems based on deep neural networks. The local parameterization of the external time-dependent inputs is defined over a set of discrete time instances and conducted using a chosen local basis over time. The resulting piecewise local parametric systems in each time interval can be learned via  our previous framework DRIPS~\cite{lu2022drips}. Once the local surrogate model is successfully constructed from training data, predictions of different initial conditions and time-dependent external inputs can be made by iterative computation of the interpolated parametric surrogate models over each discrete time instances. Our framework serves as a more data-efficient alternative to the DNN based learning method~\cite{qin2021data} as further studied in numerical experiments.

%% paper structure
The remainder of the paper is organized as follows. The problem of interest is formulated in section~\ref{sec:2}. Then we present the description of proposed methodology in section~\ref{sec:3}, which consists of the transformation induced by local parameterization and an overview of DRIPS for learning the modified system. In section~\ref{sec:4}, we test our framework on the numerical examples presented in~\cite{qin2021data}, which includes linear and nonlinear dynamical systems as well as a PDE.  The comparison with DNN demonstrates the efficiency and robustness of our method.

\section{Problem Formulation}
\label{sec:2}
We consider a multi-physics system described by coupled partial differential equations (PDEs)
\begin{align}\label{eq:genPDE}
\begin{aligned}
& \frac{\partial s_i}{\partial t} = \phi_i(\mathbf s; \gamma(\mathbf x, t)), \quad (\mathbf x,t) \in \mathcal D \times (0,T], \\
&s_i(0) = s_i^0,\qquad i = 1,\dots,N_s,
\end{aligned}
\end{align}
where the $N_s$ state variables $\mathbf s(\mathbf x,t, \gamma(\mathbf x, t)) = \{s_1,\dots,s_{N_q} \}$ vary in space $\mathbf x \in \mathcal D$ and time $t \in [0,T]$ throughout the simulation domain $\mathcal D$ during the simulation time interval $[0,T]$,  $\phi_i$ are (linear/nonlinear) differential operators that contain spatial derivatives, and $\gamma(\mathbf x,t)$ is known time-dependent inputs. When solved numerically, the spatial domain $\mathcal D$ is discretized into $N_\text{el}$ elements/nodes, leading to the discretized state variable $\mathbf S(t;\Gamma(t))$ of (high) dimension $N_S = N_s\times N_\text{el}$, satisfying the following nonautonomous ordinary differential equations (ODEs)~\eqref{eq:ODEs}

\begin{equation}\label{eq:ODEs}
\left\{
\begin{aligned}
&\frac{\text d \mathbf S}{\text dt}(t) = \boldsymbol \Phi (\mathbf S, \Gamma(t)),\\
&\mathbf S(0) = \mathbf S_0,
\end{aligned}
\right.
\end{equation}
where $\Gamma(t)$ is a vector-value time-dependent input known from spatial discretization of $\gamma(\mathbf x,t)$.

Given a set of measurement data of the state variables (elaborated in section~\ref{sec:dataset}), our goal is to construct a numerical surrogate model which can learn the behavior of the unknown nonautonomous system~\eqref{eq:ODEs}. Specifically, an accurate prediction $\hat{\mathbf S}$ of the true state variable $\mathbf S$ (satisfying~\eqref{eq:ODEs}) can be made for an arbitrary initial condition $\mathbf S_0$ and an external input process $\Gamma(t)$ within a finite time horizon $T>0$, i.e.,
\begin{equation}
\hat{\mathbf S}(t_k;\mathbf S_0,\Gamma(t))\approx \mathbf S(t_k;\mathbf S_0,\Gamma(t)),\qquad k = 1,\cdots,N_T,\quad 0 = t_0<\cdots < t_{N_T} = T.
\end{equation}
In the following context, we take $\Gamma(t)$ as a scalar function for illustration purpose. The method can be easily applied to vector-valued time-dependent inputs component by component.

\section{Methodology}
\label{sec:3}
In this section, we present the description of a framework for learning nonautonomous systems~\eqref{eq:ODEs} using DMD. It includes two major steps: 1. decomposing the dynamical system into a modified system comprising a sequence of local systems by parameterizing the external input $\Gamma(t)$ locally in time; 2. learning the local parametric systems via DRIPS. The first step follows the same reconstruction as in~\cite{qin2021data} and the second step is adopted from the DRIPS framework proposed in our previous work~\cite{lu2022drips}.

\subsection{Local parameterization and modified system}
The corresponding discrete-time dynamical system of~\eqref{eq:ODEs} is described by
\begin{equation}\label{eq:HFM-dis}
\begin{aligned}
\mathbf S(t_{k+1}) &= \boldsymbol \Phi_{\Delta t}(\mathbf S(t_{k}),\Gamma(t_k))\\
& := \mathbf S(t_{k})+\int_{t_k}^{t_k+\Delta t} \boldsymbol \Phi(\mathbf S(\tau),\Gamma(\tau)) \text d\tau,\\
&=\mathbf S(t_{k})+\int_{0}^{\Delta t} \boldsymbol \Phi(\mathbf S(t_k+\tau),\Gamma(t_k+\tau)) \text d\tau,
\end{aligned}
\end{equation}
for the uniform time discretization $t_k = k\Delta t\in [0,T]$ with $k = 0,\dots, N_T$. In each time interval $[t_k,t_{k+1}], k = 0,\cdots, N_T-1$, a local parameterization is made for the given external input function $\Gamma(t)$ in the form
\begin{equation}\label{eq:local_par}
\tilde \Gamma_k(\tau; \mathbf p_k):= \sum_{j=1}^{N_\text{par}}p_k^j b_j(\tau)\approx \Gamma(t_k+\tau),\qquad \tau \in [0,\Delta t],
\end{equation}
where $\{b_j(\tau),j = 1,\cdots, N_\text{par}\}$ is a set of prescribed analytical basis functions and 
\begin{equation}\label{eq:par}
\mathbf p_k = (p_k^1,\cdots, p_k^{N_\text{par}})\in \mathbb R^{N_\text{par}}
\end{equation}
are the basis coefficients parameterizing the local input $\Gamma(t)$ in $[t_k,t_{k+1}]$. Examples of local parameterization of a given input $\Gamma(t)$ include interpolating polynomials, Taylor polynomials and etc. (see section 3.1 in~\cite{qin2021data}). Then a global parameterization can be constructed for the external input $\Gamma(t)$ as follows
\begin{equation}\label{eq:global_par}
\tilde \Gamma(t;\mathbf p) = \sum_{k = 0}^{N_T-1}\tilde \Gamma_k(t-t_k;\mathbf p_k)\mathbb I_{[t_k,t_{k+1}]}(t),
\end{equation}
where
\begin{equation}
\mathbf p = \{\mathbf p_k\}_{k=0}^{N_T-1}\in \mathbb R^{N_T\times N_\text{par}}
\end{equation}
is a global parameter set for $\tilde \Gamma (t)$ and $\mathbb I_{[a,b]}$ is the indicator function 
\begin{equation}
\mathbb I_{[a,b]} (t) = \left\{
\begin{aligned}
&1&&\text{if }t\in [a,b],\\
&0&&\text{otherwise}.
\end{aligned}
\right.
\end{equation}

A modified system corresponding to the true (unknown) system~\eqref{eq:ODEs} is defined as follows:
\begin{equation}\label{eq:NODE_modified}
\left\{
\begin{aligned}
&\frac{\text d \tilde{\mathbf S}}{\text dt}(t) =  \boldsymbol \Phi (\tilde{\mathbf S}, \tilde\Gamma(t;\mathbf p)),\\
&\tilde{\mathbf S}(0) = \mathbf S_0,
\end{aligned}
\right.
\end{equation}
where $\tilde\Gamma(t;\mathbf p)$ is the globally parameterized input defined in~\eqref{eq:global_par}. Note that when the system input $\Gamma (t)$ is already known or given in a parametric form, i.e., $\tilde \Gamma (t) = \Gamma (t)$, the modified system~\eqref{eq:NODE_modified} is equivalent to the original system~\eqref{eq:ODEs}. When the parameterized process $\tilde \Gamma (t)$ needs to be numerically constructed, the modified system~\eqref{eq:NODE_modified} becomes an approximation to the true system~\eqref{eq:ODEs}. The approximation accuracy depends on the accuracy in $\tilde \Gamma(t)\approx \Gamma(t)$.

Following Lemma 3.1 in~\cite{qin2021data}, one can show that there exists a function $\tilde { \boldsymbol \Phi}_{\Delta t} : \mathbb R^{N_S}\times \mathbb R^{N_\text{par}} \to \mathbb R^{N_S}$, which depends on $ \boldsymbol \Phi$, such that, for any time interval $[t_k, t_{k+1}]$, the solution of ~\eqref{eq:NODE_modified} satisfies
\begin{equation}\label{eq:modified-HFM-dis}
\tilde{\mathbf S}(t_{k+1}) = \tilde { \boldsymbol \Phi}_{\Delta t} (\tilde{\mathbf S}(t_{k}),\mathbf p_k),\qquad k = 0,\cdots, N_T-1,
\end{equation}
where $\mathbf p_k$ is the local parameter set~\eqref{eq:par} for the locally parameterized input $\tilde \Gamma_k(t)$ in~\eqref{eq:local_par}.

\subsection{Learning the modified systems via DRIPS}
The function $\tilde{ \boldsymbol \Phi}_{\Delta t}$ in~\eqref{eq:modified-HFM-dis} governs the evolution of the solution of the modified system~\eqref{eq:NODE_modified} and is the target function to learn. The difficulty in learning nonautonomous system is now shifted to the task of learning the parametric system~\eqref{eq:modified-HFM-dis} in any local time interval $[t_k,t_{k+1}]$. This task falls into the category of problems where DRIPS framework applies. 

\subsubsection{Training and testing dataset}\label{sec:dataset}
We first elaborate on how to collect the training dataset. Given a set of known external inputs $\{\Gamma^{(i)}(t)\}_{i=1}^{i=N_\text{traj}}$ and  a prescribed set of discrete time instances $0 = t_0<t_1<\cdots t_k <\cdots <t_{N_{T^{(i)}}} = T^{(i)}$ with $\Delta t = t_{k+1}-t_k, k = 0,\cdots, N_{T^{(i)}}-1$, we assume $N_\text{snap}^{(k,i)}$ pairs of the input ($\mathbf S_j^{(i)}(t_k), \Gamma^{(i)}(t_k)$)-output ($\mathbf S_j^{(i)}(t_{k+1})$) responses from the true discrete-time dynamical system~\eqref{eq:HFM-dis} are available along the $i$th trajectory within the time interval $[t_k,t_{k+1}]$, i.e.,
\begin{equation}
\mathbf S_j^{(i)}(t_{k+1}) = \boldsymbol \Phi_{\Delta t}(\mathbf S_j^{(i)}(t_k), \Gamma^{(i)}(t_k)),\qquad j = 1,\cdots, N_\text{snap}^{(k,i)},\quad k = 0,\cdots N_{T^{(i)}}-1, \quad i = 1,\cdots, N_\text{traj}.
\end{equation}

The local parameterization of $\Gamma^{(i)}(t_k)$ gives $\tilde \Gamma_k^{(i)}(\tau;\mathbf p_k^{(i)})$, where $\tau\in [0, \Delta t]$ and $\mathbf p_k^{(i)}$ is the parameter set for the local parameterization of the input in the form of~\eqref{eq:local_par}. Now denote each local dataset as
\begin{equation}
\mathcal S_\text{train}^{(k,i)} = \{(\mathbf S^{(i)}_j(t_k),\mathbf p_k^{(i)}),\mathbf S^{(i)}_j(t_{k+1})\},\quad j = 1,\cdots, N^{(k,i)}_\text{snap}.
\end{equation}

Notice that the following input ($\mathbf S^{(i)}_j(t_k),\mathbf p_k^{(i)}$)-output $\mathbf S^{(i)}_j(t_{k+1})$ responses hold approximately in the underlying modified system~\eqref{eq:modified-HFM-dis}, i.e.,
\begin{equation}\label{eq:train-HFM}
\mathbf S^{(i)}_j(t_{k+1}) \approx \tilde { \boldsymbol \Phi}_{\Delta t} (\mathbf S^{(i)}_j(t_{k}),\mathbf p^{(i)}_k),\qquad j = 1,\cdots, N^{(k,i)}_\text{snap},\quad k = 0,\cdots N_{T^{(i)}}-1, \quad i = 1,\cdots, N_\text{traj}.
\end{equation}

Then the full training dataset can be represented as 
\begin{equation}
\mathcal S_\text{train} = \cup_{k,i} \mathcal S_\text{train}^{(k,i)}, \quad k = 0,\cdots N_{T^{(i)}}-1, \quad i = 1,\cdots, N_\text{traj}.
\end{equation}

On re-ordering using a single index, the training dataset takes the form
\begin{equation}\label{eq:train_data}
\mathcal S_\text{train} = \cup_{m=1}^{N_\text{MC}} \mathcal S_\text{train}^{(m)}, \quad N_\text{MC} =\sum_{i=1}^{N_\text{traj}} N_{T^{(i)}},
\end{equation}
where 
\begin{equation}\label{eq:train_data_sub}
\mathcal S_\text{train}^{(m)} = \{(\mathbf S^{(m)}_j(0),\mathbf p^{(m)}),\mathbf S^{(m)}_j(\Delta t)\},\quad j = 1,\cdots, N^{(m)}_\text{snap}.
\end{equation}

Our goal is to construct from this dataset $\mathcal S_\text{train}$ a surrogate model of the unknown system~\eqref{eq:modified-HFM-dis}. This will allow us to predict the trajectory of the state variable $\mathbf S(t;\Gamma^*(t))$ from a different external input $\Gamma^*(t)\notin \{\Gamma^{(i)}(t)\}_{i=1}^{i=N_\text{traj}}$, i.e., $\{\mathbf S(t_0,\Gamma^*(t_0)),\cdots, \mathbf S(t_{N_{T^*}},\Gamma^*(t_{N_{T^*}}))\}$ at a low cost. The corresponding test dataset can be obtained after representing $\Gamma^*(t_k)$ by $\tilde\Gamma_k^*(\tau;\mathbf p_k^*)$ via local parameterization~\eqref{eq:local_par},
\begin{equation}\label{eq:data_test}
\mathcal S_\text{test} = \{\tilde{\mathbf S}(t_0,\mathbf p_0^*),\cdots, \tilde{\mathbf S}(t_{N_{T^*}},\mathbf p_{N_{T^*}}^*)\}.
\end{equation}

The target function $\tilde{ \boldsymbol \Phi}_{\Delta t}$ and the modified system~\eqref{eq:modified-HFM-dis} can be learned within our DRIPS framework, which consists of the offline step and online step that are detailed in Sections~\ref{sec:offline} and~\ref{sec:online}, respectively.

\subsubsection{Offline Step: DMD-Based Surrogates}
\label{sec:offline}

For nonautonomous system~\eqref{eq:ODEs}, the dynamics of $\mathbf S$, i.e., the functional form of $ \boldsymbol \Phi$ in~\eqref{eq:ODEs} or $ \boldsymbol \Phi_{\Delta t}$ in~\eqref{eq:HFM-dis} is unknown. So is $\tilde{ \boldsymbol \Phi}_{\Delta t}$ in~\eqref{eq:modified-HFM-dis} for the modified systems. In the first step of our algorithm, we replace the unknown discrete system~\eqref{eq:modified-HFM-dis} with its linear surrogate  model constructed  from the dataset~\eqref{eq:train_data}.  The latter task is facilitated by the Koopman operator theory, which allows one to handle the potential nonlinearity in the unknown dynamic $ \boldsymbol \Phi$ and $\tilde{ \boldsymbol \Phi}_{\Delta t}$:
\begin{definition}[Koopman operator~\cite{kutz2016dynamic}]
For nonlinear dynamic system~\eqref{eq:ODEs}, the Koopman operator $\mathcal K^{\Gamma(t)}$ is an infinite-dimensional linear operator that acts on all observable functions $g: \mathcal M\to \mathbb C$ so that
\begin{equation}\label{eq:Koopman}
\mathcal K^{\Gamma(t)} g(\mathbf S(t)) = g( \boldsymbol \Phi(\mathbf S(t),\Gamma(t))).
\end{equation}
For discrete dynamic system~\eqref{eq:modified-HFM-dis}, the discrete-time Koopman operator $\mathcal K_{\Delta t}^{\mathbf p_k}$ is
\begin{equation}\label{eq:Koopman-dis}
\mathcal K_{\Delta t} ^{\mathbf p_k}g(\mathbf S(t_k;\mathbf p_k)) = g(\tilde{ \boldsymbol \Phi}_{\Delta t}(\mathbf S(t_k;\mathbf p_k),\mathbf p_k)) = g(\mathbf S(t_{k+1};\mathbf p_k)).
\end{equation}
\end{definition}

The Koopman operator transforms the finite-dimensional nonlinear problem~\eqref{eq:modified-HFM-dis} in the state space into the infinite-dimensional linear problem~\eqref{eq:Koopman-dis} in the observable space. Since $\mathcal K_{\Delta t}^{\mathbf p_k}$ is an infinite-dimensional linear operator, it is equipped with infinite eigenvalues $\{\lambda_i(\mathbf p_k)\}_{i=1}^\infty$ and eigenfunctions $\{\phi_i(\mathbf p_k)\}_{i=1}^\infty$. In practice, one has to make do with a finite number of the eigenvalues and eigenfunctions. The following assumption is essential to both a finite-dimensional approximation and the choice of observables.

\begin{assump}\label{assump1}
Let $\mathbf y (t_k;\mathbf p_k)$ denote an $N\times 1$ vector of observables,
\begin{equation}
\mathbf y(t_k;\mathbf p_k) = \mathbf g(\mathbf S(t_k;\mathbf p_k) )= \begin{bmatrix}
g_1(\mathbf S(t_k;\mathbf p_k))\\
\vdots\\
g_N(\mathbf S(t_k;\mathbf p_k))
\end{bmatrix},
\end{equation}
where $g_j:\mathcal M\to \mathbb C$ is an observable function with $j = 1,\dots,N$. If the chosen observables $\mathbf g$ are restricted to an invariant subspace spanned by eigenfunctions of the Koopman operator $\mathcal K_{\Delta t}^{\mathbf p_k}$, then they induce a linear operator $\mathbf K(\mathbf p_k)$ that is finite-dimensional and advances these eigen-observable functions on this subspace~\cite{brunton2016koopman}.
\end{assump}
Assumption~\ref{assump1} enables one to deploy the DMD Algorithm~\ref{alg:DMD} to approximate the $N$-dimensional linear operator $\mathbf K(\mathbf p_k)$ and its low-dimensional approximation $\mathbf K_r(\mathbf p_k)$ of rank $r$. At each parameter point $\mathbf p_k^{(i)}$ with $k = 1,\cdots, N_{T^{(i)}}, i = 1,\dots, N_{\text{traj}}$, the discrete system~\eqref{eq:train-HFM} on state space is approximated by an $N$-dimensional linear surrogate model
\begin{equation}\label{eq:surrogate}
\mathbf y(t_{k+1};\mathbf p_k^{(i)}) = \mathbf K(\mathbf p_k^{(i)})\mathbf y(t_{k};\mathbf p_k^{(i)})
\end{equation}
on observable space. The two spaces are connected  by the observable function $\mathbf g$ and its inverse $\mathbf g^{-1}$. Algorithm~\ref{alg:DMD} directly induces the  ROM for~\eqref{eq:surrogate},
\begin{equation}
\label{eq:ROM}
\mathbf y_r(t_{k+1};\mathbf p_k^{(i)}) = \mathbf K_r(\mathbf p_k^{(i)})\mathbf y_r(t_{k};\mathbf p_k^{(i)}).
\end{equation}
Here $\mathbf y_r(t_{k};\mathbf p_k^{(i)})$ is the reduced-order observable vector or dimension $r$. In terms of a ROB $\mathbf V(\mathbf p_k^{(i)})$, these are expressed as
\begin{equation}
\mathbf y(t_{k};\mathbf p_k^{(i)}) = \mathbf V (\mathbf p_k^{(i)})\mathbf y_r(t_{k};\mathbf p_k^{(i)})
\quad\text{and}\quad 
\mathbf K_r(\mathbf p_k^{(i)}) =  \mathbf V (\mathbf p_k^{(i)}) ^\top \mathbf K(\mathbf p_k^{(i)})\mathbf V (\mathbf p_k^{(i)}).
\end{equation}

The following Algorithm~\ref{alg:DMD} can be trained on the re-ordered dataset~\eqref{eq:train_data} then,

\begin{algorithm}[h]
Input: $\mathcal S_\text{train}^{(m)}$ in~\eqref{eq:train_data_sub}, observable function $\mathbf g$
\begin{enumerate}
\item Create data matrices of the observables
\begin{equation}
\mathbf Y_1(\mathbf p^{(m)}) = \begin{bmatrix}
|&&|\\
\mathbf g(\mathbf S_1(0;\mathbf p^{(m)}))&\dots&\mathbf g(\mathbf S_{N^{(m)}_\text{snap}}(0;\mathbf p^{(m)}))\\
|&&|
\end{bmatrix}, 
\end{equation}

\begin{equation}
\mathbf Y_2(\mathbf p^{(m)}) = \begin{bmatrix}
|&&|\\
\mathbf g(\mathbf S_1(\Delta t;\mathbf p^{(m)}))&\dots&\mathbf g(\mathbf S_{N^{(m)}_\text{snap}}(\Delta t;\mathbf p^{(m)}))\\
|&&|
\end{bmatrix}.
\end{equation}

\item Apply SVD $\mathbf Y_1(\mathbf p^{(m)}) \approx \mathbf V(\mathbf p^{(m)}) \boldsymbol \Sigma(\mathbf p^{(m)}) \mathbf Z (\mathbf p^{(m)}) ^*$, where  $\mathbf V(\mathbf p^{(m)}) \in \mathbb C^{N\times r}$, $\boldsymbol \Sigma(\mathbf p^{(m)})  \in \mathbb C^{r\times r}$, $\mathbf Z(\mathbf p^{(m)}) \in \mathbb C^{r\times N^{(m)}_\text{snap}}$, and $r$ is the truncation rank chosen by certain criteria and kept the same for all $m= 1,\dots, N_\text{MC}$.
\item Compute $\mathbf K_r(\mathbf p^{(m)}) = \mathbf V(\mathbf p^{(m)})^*\mathbf Y_2(\mathbf p^{(m)})\mathbf Z(\mathbf p^{(m)})\boldsymbol \Sigma(\mathbf p^{(m)})^{-1}$ as an $r\times r$ low-rank approximation of $\mathbf K(\mathbf p^{(m)})$.
\item Compute $\mathbf P^{(m_1,m_2)} = \mathbf V(\mathbf p^{(m_1)})^\top\mathbf V(\mathbf p^{(m_2)})$ for $m_1,m_2 = 1,\dots, N_\text{MC}$.
\end{enumerate}
Output: $\mathbf V (\mathbf p^{(m)})$, $\mathbf K_r(\mathbf p^{(m)})$ and $\mathbf P^{(m_1,m_2)}$.
\caption{DMD algorithm on observable space~\cite{kutz2016dynamic} for parameter point $\mathbf p^{(m)}$, $ m = 1,\dots, N_{\text{MC}}$.}
\label{alg:DMD}
\end{algorithm}

\begin{remark}
The construction of DMD surrogates can be done offline, which allows one to precompute $\mathbf P^{(m_1,m_2)}$ for the later online step. Although this offline step takes a majority of the computational time in the whole framework, mostly due to the computation of the high-fidelity training data~\eqref{eq:train_data}, its output can be precomputed and stored efficiently; the output storage is $(N \cdot r + r \cdot r + r \cdot r \cdot (N_\text{MC}+1)/2) \cdot N_\text{MC}$. %In practice, real-time applications can be enabled as long as the online step is sufficient enough.
\end{remark}

\begin{remark}
A theorem in~\cite{tu2013dynamic} establishes connections between the DMD theory and the Koopman spectral analysis under specific conditions on the observables and collected data. This theorem indicates that judicious selection of the observables is critical to the success of the Koopman method. There is no principled way to select observables without expert knowledge of a dynamical system. Machine learning techniques can be deployed to identify relevant terms in the dynamics from data, which then guides the selection of the observables~\cite{schmidt2009distilling,wang2011predicting}. In our numerical examples, we rely on knowledge of the underlying physics to select the observables, as was done, e.g., in~\cite{lu2020lagrangian,lu2020prediction,lu2021dynamic,lu2021extended,%williams2014kernel,
li2017extended}).
\end{remark}

\subsubsection{Online Step: Interpolation of ROBs and PROMs}
\label{sec:online}

For a different external input $\Gamma^*(t)\notin \{\Gamma^{(i)}(t)\}_{i=1}^{N_\text{traj}}$, the goal is to compute $\{\mathbf S(t_0,\Gamma^*(t_0)),\cdots, \mathbf S(t_{N_{T^*}},\Gamma^*(t_{N_{T^*}}))\}$, at a low cost without evaluating~\eqref{eq:ODEs}. Using the same local parameterization as in~\eqref{eq:local_par}, we seek to approximate the test dataset~\eqref{eq:data_test} via the PROM
\begin{equation}\label{eq:ROM-iter}
\mathbf y_r(t_{k+1};\mathbf p_k^{*}) = \mathbf K_r(\mathbf p_k^{*})\mathbf y_r(t_{k};\mathbf p_k^{*}).
\end{equation}
Subsequently, the observables $\mathbf y$ and the state variable $\mathbf S$ are estimated as
\begin{equation}\label{eq:ROM-state}
\mathbf y(t_{k};\mathbf p_k^*) = \mathbf V (\mathbf p_k^*)\mathbf y_r(t_{k};\mathbf p_k^*),\qquad
\mathbf S(t_k; \mathbf p_k^*) = \mathbf g^{-1}(\mathbf y(t_k;\mathbf p_k^*)).
\end{equation}
 Therefore, the online task comprises the computation of three quantities, $\mathbf V (\mathbf p_k^*)$, $\mathbf K_r (\mathbf p_k^*)$, and $\mathbf y_r(t_k;\mathbf p_k^*)$. Notice that in general $\mathbf p_k^*\notin \{\mathbf p^{(m)}\}_{m=1}^{N_\text{MC}}$ for $k = 0,\cdots, N_{T^*}$.

\paragraph{ROB Interpolation}
\label{sec:interpROB}

We rely on interpolation on the Grassman manifold to compute the ROB $\mathbf V (\mathbf p_k^*)$ from the dataset $\{ \mathbf V (\mathbf p^{(1)}),\dots,\mathbf V (\mathbf p^{(N_\text{MC})})\}$. We briefly review this interpolation approach~\cite{amsallem2008interpolation} below.

\begin{definition} The following manifolds are of interest:
\begin{itemize}
\item Grassmann manifold $\mathcal G(r,N)$ is the set of all subspaces in $\mathbb R^N$  of dimension $r$;
\item Orthogonal Stiefel manifold $\mathcal S\mathcal T(r,N)$ is the set of orthogonal ROB matrices in $\mathbb R^{r\times N}$. 
\end{itemize}
\end{definition}

The ROB $\mathbf V (\mathbf p^{(m)}) \in \mathbb R^{N\times r}$, with $m = 1,\dots, N_\text{MC}$ and $r\leq N$, is the full-rank column matrix, whose columns provide a basis of subspace $\mathcal S_m$ of dimension $r$ in $\mathbb R^N$. While an associated ROM is not uniquely defined by the ROB, it is uniquely defined by the subspace $\mathcal S_m$. Therefore, the correct entity to interpolate is the subspace $\mathcal S_m$, rather than the ROB  $\mathbf V (\mathbf p^{(m)})$. Hence, the goal is to compute $\mathcal S_{k,*} = \text{range}(\mathbf V(\mathbf p_k^*))$ by interpolating between $\{\mathcal S_m\}_{m=1}^{N_{\text{MC}}}$ with access to the ROB $\mathbf V(\mathbf p_k^*)$. 

The subspaces $\mathcal S_1,\dots,S_{N_\text{MC}}$ belong to the Grassmann manifold $\mathcal G(r,N)$~\cite{absil2004riemannian,boothby2003introduction,helgason2001differential,rahman2005multiscale,edelman1998geometry} . Each $r$-dimensional subspace $\tilde {\mathcal S}$ of $\mathbb R^N$ is regarded as a point of $\mathcal G(r,N)$ and is nonuniquely represented by a matrix $\tilde{\mathbf V}\in \mathbb R^{N\times r}$, whose columns span the subspace $\tilde{\mathcal S}$. The matrix $\tilde{\mathbf V}$ is chosen among those whose columns form a set of orthonormal vectors in $\mathbb R^N$ and belong to the orthogonal Stiefel manifold $\mathcal S\mathcal T(r,N)$~\cite{absil2004riemannian,edelman1998geometry}. There exists a projection map~\cite{absil2004riemannian} from $\mathcal S\mathcal T(r,N)$ onto $\mathcal G(r,N)$, as each matrix in $\mathcal S\mathcal T(r,N)$ spans an $r$-dimensional  subspace of $\mathbb R^N$ and different matrices can span the same subspace. At each point $\tilde{\mathcal S}$ of the manifold $\mathcal G(r,N)$, there exists a tangent space $\mathcal T_{\tilde{\mathcal S}}$~\cite{absil2004riemannian,edelman1998geometry} of the same dimension~\cite{edelman1998geometry}. Each point in this space is represented by a matrix $\tilde {\mathbf M}\in \mathbb R^{N\times r}$. Since $\mathcal T_{\tilde{\mathcal S}}$ is a vector space, usual interpolation is allowed for the matrices representing its points. Let $\mathbf M^{m} = m_{\mathbf V}(\mathbf V(\mathbf p^{(m)}))$, where $m_{\mathbf V}$ denotes the map from the matrix manifolds $\mathcal G(r,N)$ onto the tangent space $\mathcal T_{\tilde{\mathcal S}}$. This suggests a strategy of computing $\mathbf M^{k,*}$ via usual interpolation between $\{\mathbf M^m\}_{m=1}^{N_{\text{MC}}}$ and then evaluating $\mathbf V(\mathbf p_k^*)$ through the inverse map $m_{\mathbf V}^{-1}(\mathbf M^{k,*})$.

The map $m_{\mathbf V}$ is chosen to be the logarithmic mapping, which maps the Grassmann manifold onto its tangent space, and $m_{\mathbf V}^{-1}$ is chosen to be the exponential mapping, which maps the tangent space onto the Grassmann manifold itself. This choice borrows concepts of geodesic path on a Grassmann manifold from differential geometry~\cite{absil2004riemannian,boothby2003introduction,wald1984general,do1992riemannian}. This strategy, discussed in detail in~\cite{amsallem2008interpolation}, is implemented in Algorithm~\ref{alg:ROB}.

\begin{algorithm}[h]
Input: $\{ \mathbf V (\mathbf p^{(m)})\}_{m=1}^{N_{\text{MC}}}$, $\{ \mathbf P^{(m_1,m_2)}\}_{m_1,m_2=1}^{N_{\text{MC}}}$, $\{\mathbf p^{(m)}\}_{m=1}^{N_{\text{MC}}}$ and target parameter point $\{\mathbf p_k^*\}_{k=0}^{N_{T^*}}$
\begin{enumerate}
\item Denote $\mathcal S_m = \text{range}(\mathbf V(\mathbf p^{(m)}))$, for $m = 1,\dots N_{\text{MC}}$. A point $\mathcal S_{m_0}$ with $m_0\in \{1,\dots N_{\text{MC}}\}$ of the manifold is chosen as a reference and origin point for interpolation.
\item Select points $\mathcal S_m$ with $m \in \mathcal I_{m_0}\subset \{1,\dots, N_{\text{MC}}\}$, which lie in a sufficiently small neighborhood of $\mathcal S_{m_0}$, and use the logarithm map $\log_{\mathcal S_{m_0}}$ to map $\{\mathcal S_m\}_{m\in \mathcal I_{m_0}}$ onto matrices $\{\mathbf M^m\}_{m\in \mathcal I_{m_0}}$ representing the corresponding points of $\mathcal T_{\mathcal S_{m_0}}$. This is computed as
\begin{equation}
\begin{aligned}
&(\mathbf I -\mathbf V(\mathbf p^{(m_0)})\mathbf V(\mathbf p^{(m_0)})^\top)\mathbf V(\mathbf p^{(m)})(\mathbf P^{(m_0,m)})^{-1} = \mathbf U_m\boldsymbol \Omega_m\mathbf W_m^\top,\quad \text{(thin SVD)}\\
&\mathbf M^m = \mathbf U_m\tan^{-1}(\boldsymbol \Omega_m)\mathbf W_m^\top.
\end{aligned}
\end{equation}
\item Compute $\mathbf M^{k,*}$ by interpolating $\{\mathbf M^m\}_{m\in \mathcal I_{m_0}}$ entry by entry:
\begin{equation}\label{eq:entry-by-entry}
M_{ij}^{k,*} = \mathcal P(\mathbf p_k^*;\{M_{ij}^{m},\mathbf p^{(m)}\}_{m\in \mathcal I_{m_0}}), \quad i =1,\dots, N,\quad j = 1,\dots, r. 
\end{equation}
\item Use the exponential map $\exp_{\mathcal S_{m_0}}$ to map the matrix $\mathbf M^{k,*}$, representing a point of $\mathcal T_{\mathcal S_{m_0}}$, onto the desired subspace $\mathcal S_{k,*}$ on $\mathcal G(r,N)$ spanned by the ROB $\mathbf V(\mathbf p_k^*)$. This is computed as
\begin{equation}
\begin{aligned}
&\mathbf M^{k,*} =  \mathbf U_{k,*}\tan^{-1}(\boldsymbol \Omega_{k,*})\mathbf W_{k,*}^\top,\quad \text{(thin SVD)}\\
&\mathbf V(\mathbf p_k^*) = \mathbf V(\mathbf p^{(m_0)})\mathbf W_{k,*}\cos(\boldsymbol \Omega_{k,*})+\mathbf U_{k,*}\sin(\boldsymbol \Omega_{k,*}).
\end{aligned}
\end{equation}
\end{enumerate}
Output: $\{\mathbf V(\mathbf p_k^*)\}_{k=1}^{N_{T^*}}$
\caption{Interpolation of ROBs~\cite{amsallem2008interpolation}}
\label{alg:ROB}
\end{algorithm}

\begin{remark}
The choice of the interpolation method $\mathcal P$ depends on the  dimension of the parameter space, $N_\text{par}$. If $N_\text{par}=1$, a univariate (typically, a Lagrange type) interpolation method is chosen. Otherwise, a multivariate interpolation scheme (see, e.g.,~\cite{spath1995one,de1992computational}) is chosen.
\end{remark}

\begin{remark}
Since the logarithmic map $\log_{\mathcal S_{m_0}}$ is defined in a neighborhood of $\mathcal S_{m_0}\in \mathcal G(r,N)$, the method is expected to be insensitive to the choice of the reference point $\mathcal S_{m_0}$ in step 1 of Algorithm~\ref{alg:ROB}. This is confirmed in numerical experiments~\cite{amsallem2008interpolation}.
\end{remark}

\paragraph{PROM Interpolation}

The reduced-order operator $\mathbf K_r (\mathbf p_k^*)$ in~\eqref{eq:ROM-iter} is computed via interpolation on the matrix manifold between the ROMs $\{ \mathbf K_r (\mathbf p^{(1)}),\dots,\mathbf K_r (\mathbf p^{(N_\text{MC})})\}$. This is done in two steps~\cite{amsallem2011online}: 
\begin{itemize}
\item Step A). Since any ROM can be endowed with multiple alternative ROBs, the resulting ROMs may have been computed in different generalized coordinates system. The validity of an interpolation may crucially depend on the choice of the representative element within each equivalent class. Given the precomputed ROMs $\{ \mathbf K_r (\mathbf p^{(1)}),\dots,\mathbf K_r (\mathbf p^{(N_\text{MC})})\}$, a set of congruence transformations is determined so that a representative element of the equivalent ROBs for each precomputed ROM is chosen to assign the precomputed ROMs into consistent sets of generalized coordinates. The consistency is enforced by solving the orthogonal Procrustes problems~\cite{van1996matrix},
\begin{equation}
\min_{\mathbf S_m,\mathbf S_m^\top\mathbf S_m = \mathbf I_r}\|\mathbf V(\mathbf p^{(m)})^\top\mathbf S_m-\mathbf V(\mathbf p^{(m_0)})\|_F, \qquad
\forall m = 1,\dots, N_\text{MC},
\end{equation}
where $m_0\in \{1,\dots N_{\text{MC}}\}$ is chosen as a reference configuration. The representative element is identified by solving the above problem analytically. This procedure is summarized in Algorithm~\ref{alg:stepA}. 

\begin{algorithm}[h]
Input: $\{ \mathbf K_r (\mathbf p^{(1)}),\dots,\mathbf K_r (\mathbf p^{(N_\text{MC})})\}$, $\{ \mathbf P^{(i,j)}\}_{m_1,m_2=1}^{N_{\text{MC}}}$, reference configuration choice $m_0$

\textbf{For}  $m\in \{1,\dots N_{\text{MC}}\}\setminus\{m_0\}$
\begin{itemize}
\item Compute $\mathbf P^{(m,m_0)} = \mathbf U_{m,m_0}\boldsymbol \Sigma_{m,m_0}\mathbf Z_{m,m_0}^\top$ (SVD),
\item Compute $\mathbf S_m = \mathbf U_{m,m_0}\mathbf Z_{m,m_0}^\top$,
\item Transform $\tilde{\mathbf K}_r (\mathbf p^{(m)}) = \mathbf S_m^\top \mathbf K_r (\mathbf p^{(m)}) \mathbf S_m$
\end{itemize}
\textbf{End}

Output: $\{ \tilde{\mathbf K}_r (\mathbf p^{(1)}),\dots,\tilde{\mathbf K}_r (\mathbf p^{(N_\text{MC})})\}$
\caption{Step A of the PROM interpolation~\cite{amsallem2008interpolation}}
\label{alg:stepA}
\end{algorithm}

\begin{remark}
An optimal choice of the reference configuration $m_0$, if it exists, remains an open problem.
\end{remark}
\begin{remark}
Step A is related to mode-tracking procedures based on the modal assurance criterion (MAC)~\cite{ewins2009modal}. This connection is illucidated in~\cite{amsallem2011online}.
\end{remark}

\item Step B). The transformed ROMs $\{ \tilde{\mathbf K}_r (\mathbf p^{(1)}),\dots,\tilde{\mathbf K}_r (\mathbf p^{(N_\text{MC})})\}$ are interpolated to compute ROMs $\{\mathbf K_r (\mathbf p_k^*)\}_{k=1}^{N_{T^*}}$. Similar to Section~\ref{sec:interpROB}, this interpolation must be performed on a specific manifold  containing both $\{ \tilde{\mathbf K}_r (\mathbf p^{(1)}),\dots,\tilde{\mathbf K}_r (\mathbf p^{(N_\text{MC})})\}$ and $\{\mathbf K_r (\mathbf p_k^*)\}_{k=1}^{N_{T^*}}$, so that the distinctive properties (e.g., orthogonality, nonsingularity) are preserved. The main idea again is to first map all the precomputed matrices onto the tangent space to the matrix manifold of interest at a chosen reference point using the logarithm mapping, then interpolate the mapped data in this linear vector space, and finally map the interpolated result back onto the manifold of interest using the associated exponential map. This is done in Algorithm~\ref{alg:stepB}.

\begin{algorithm}[h]
Input: $\{ \tilde{\mathbf K}_r (\mathbf p^{(1)}),\dots,\tilde{\mathbf K}_r (\mathbf p^{(N_\text{MC})})\}$, reference configuration choice $m_0$
\begin{enumerate}
\item \textbf{For}  $m\in \{1,\dots N_{\text{MC}}\}\setminus\{m_0\}$
\begin{itemize}
\item Compute $\mathbf M^m = \log_{\tilde{\mathbf K}_r(\mathbf p^{(m_0)})}(\tilde{\mathbf K}_r(\mathbf p^{(m)}))$
\end{itemize}
\textbf{End}
\item Compute $\mathbf M^{k,*}$ by interpolating $\{\mathbf M^m\}_{m\in \mathcal I_{m_0}}$ entry by entry, as in~\eqref{eq:entry-by-entry}
\item Compute $\mathbf K_r(\mathbf p_k^*) = \exp_{\tilde{\mathbf K}_r(\mathbf p^{(m_0)})}(\mathbf M^{k,*})$
\end{enumerate}
Output: $\{\mathbf K_r (\mathbf p_k^*)\}_{k=1}^{N_{T^*}}$
\caption{Step B of the PROM interpolation~\cite{amsallem2008interpolation}}
\label{alg:stepB}
\end{algorithm}

\begin{remark}
The $\log$ and $\exp$ in Algorithm~\ref{alg:stepB} denote the matrix logarithm and exponential respectively. The specific expressions of different matrix manifolds of interest are listed in Table 4.1 of~\cite{amsallem2011online}.
\end{remark}
\end{itemize}

\paragraph{Computation of the Solution}
With $\{\mathbf K_r (\mathbf p_k^*)\}_{k=0}^{N_{T^*}}$ and $\{\mathbf V(\mathbf p_k^*)\}_{k=0}^{N_{T^*}}$ computed, we use~\eqref{eq:ROM-iter} and~\eqref{eq:ROM-state} to obtain the solution.

\subsubsection{Algorithm Summary}
The proposed framework is implemented in Algorithm~\ref{alg:whole}.

\begin{algorithm}[H]
\textit{Offline Step:} 

\textbf{For} $m = 1,\dots,N_\text{MC}$,
$$
\begin{aligned}
&\text{Get the training data~\eqref{eq:train_data}, } \\
&\text{Input:} \ \mathcal S_\text{train}^{(m)} \ \text{and} \ \mathbf g\xrightarrow{\text{Algorithm~\ref{alg:DMD}}} \text{Output:} \ \mathbf V(\mathbf p^{(m)}),  \mathbf K_r(\mathbf p^{(m)}) \ \text{and} \ \mathbf P^{(m_1,m_2)}
\end{aligned}$$

\textbf{End}

\textit{Online Step:}
\begin{itemize}
\item Interpolation of ROBs:
$$\text{Input:} \ \{ \mathbf V (\mathbf p^{(m)})\}_{m=1}^{N_{\text{MC}}}, \ \{ \mathbf P^{(m_1,m_2)}\}_{m_1,m_2=1}^{N_{\text{MC}}}, \ \{\mathbf p^{(m)}\}_{m=1}^{N_{\text{MC}}}, \  \{\mathbf p_k^*\}_{k=0}^{N_{T^*}}\xrightarrow{\text{Algorithm~\ref{alg:ROB}}} \text{Output:} \ \{\mathbf V(\mathbf p_k^*)\}_{k=1}^{N_{T^*}}$$
\item Interpolation of PROMs:
$$\text{Input:} \ \{ \mathbf K_r (\mathbf p^{(m)})\}_{m=1}^{N_\text{MC}}, \ \{ \mathbf P^{(m_1,m_2)}\}_{m_1,m_2=1}^{N_{\text{MC}}}, \  \text{reference choice} \ m_0\xrightarrow{\text{Algorithms~\ref{alg:stepA} \&~\ref{alg:stepB}}} \text{Output:} \ \mathbf K_r(\mathbf p_k^*)$$
\item DMD reconstruction:
$$\text{Input:} \ \mathbf K_r(\mathbf p_k^*), \ \mathbf V(\mathbf p_k^*), \ \mathbf y(t=0; \mathbf p_0^*), \ \mathbf g\xrightarrow{\eqref{eq:ROM-iter}, \eqref{eq:ROM-state}} \text{Output:} \ \mathbf S_\text{DMD}(t_k;\mathbf p_k^*)$$
\end{itemize}
\caption{Learning nonautonomous system via DMD}
\label{alg:whole}
\end{algorithm}

\begin{remark}
The sampling strategy for $\{\mathbf p^{(1)},\dots, \mathbf p^{(N_\text{MC})}\}$ in the parameter space plays a key role in the accuracy of the subspace approximation. The so-called ``curse of dimensionality", i.e., the number of training samples $N_\text{MC}$ needed grows exponentially with the number of  parameters, $N_\text{par}$, is a well-known challenge. In general, uniform sampling is used for $N_\text{par}\leq 5$ and moderately computationally intensive HFMs, latin hypercube sampling is used for $N_\text{par}>5$ and moderately computationally intensive HFMs, and adaptive, goal-oriented, greedy sampling is used for $N_\text{par}>5$ and highly computationally intensive HFMs. We limit our numerical experiments to $N_\text{par} =3$ for simplicity, leaving the  challenge posed by the curse of dimensionality for future studies.
\end{remark}

\section{Numerical Experiments}
\label{sec:4}
In this section, we tested our proposed framework on every numerical example presented in~\cite{qin2021data}. Synthetic data generated from known dynamical systems with known time-dependent inputs is employed to validate the proposed framework. The training data are generated by solving the known system with a high-resolution numerical scheme in Matlab, e.g., \textit{ode45}, which is based on an explicit fourth-order Runge-Kutta formula. In all the following examples, the inputs are parameterized locally the same way as in~\cite{qin2021data}, i.e., by  interpolating polynomials over equally spaced points. Other parameterization strategies produce similar results and thus not included here for brevity.

For the convenience, we generate the training dataset in the form of~\eqref{eq:train_data} with $N_T^{(i)} \equiv 1$ for $i = 1,\cdots, N_\text{traj}$, i.e., each trajectory only contains two data points and $N_\text{MC} = N_\text{traj}$. For each subset $\mathcal S_\text{train}^{m}$, the first data entry is randomly sampled from a domain $I_{\mathbf S}$ using uniform distribution. The second data entry is produced by solving the underlying reference dynamical system with a time step $\Delta t$ and subject to a parameterized external input in the form of~\eqref{eq:local_par}. We take $\Delta t = 0.1$ in all numerical examples except when other values are specified. For simplicity, we take the same number of snapshots/data pairs in each parameter point, i.e., $N_\text{snap}^{(m)} = N_\text{snap}$ for all $m$. Instead of uniformly sampling the parameters~\eqref{eq:par} as in~\cite{qin2021data}, we take the parameter values on the Cartesian grid of the $N_\text{par}$-dimensional parameter space with $3$ points (two endpoints and one middle point) in each dimension, i.e., $N_\text{MC} =  N_\text{traj} = 3^{N_\text{par}}$. Therefore, the total amount of training data is $N_\text{snap}\times  3^{N_\text{par}}$ pairs in the form of~\eqref{eq:train_data_sub}. This number will be specified in each numerical example and we see that it is much smaller than what is needed in training a neural network ($O(10^5)$ in~\cite{qin2021data}).

The sampling domain $I_{\mathbf S}$ and parameter domain $I_{\mathbf p}$ are determined by prior knowledge of the underlying unknown equations to ensure the range of the target trajectories are covered by $I_{\mathbf S}$ for parameters in $I_{\mathbf p}$ so that the assumptions in the theoretical analysis of~\cite{qin2021data} are satisfied. The same $I_{\mathbf S}$ and $I_{\mathbf p}$ are adopted from~\cite{qin2021data} and are specified separately for each example.

Our proposed learning algorithm~\ref{alg:whole} is then applied to the training dataset. The observable $\mathbf g$ is determined by the underlying system and will be specified in each numerical example. Notice that the choice of observables may not be unique and optimal. Therefore, we consider the optimal construction of observables out of the scope of this paper and employ several standard tricks in DMD studies to construct reasonable observables. Once the learned model is constructed, we conduct system prediction iteratively using the model with new initial conditions and new external inputs. The prediction results are then compared with the reference solution obtained by solving the exact system with the same new inputs.

\subsection{Linear scalar equation with source}
We first consider the following scalar equation
\begin{equation}\label{eq:4-1}
\frac{dS}{dt} = -\alpha(t)S +\beta(t),
\end{equation}
where the time-dependent inputs $\alpha(t)$ and $\beta(t)$ are locally parameterized with polynomials of degree $2$. As a result, the local parameter set~\eqref{eq:par} $\mathbf p_k \in\mathbb R^{N_{\text{par}}}$ with $N_\text{par} = 3+3 = 6$. The training data are generated from $N_\text{MC} = N_\text{traj} = 3^6$ parameter points on the Cartesian grid of the parameter space $I_{\mathbf p} = [-5,5]^6$. On each parameter point $\mathbf p^{(m)}$, $N_\text{snap}^{(m)} = N_\text{snap} = 2$ pairs of data in the form of~\eqref{eq:train_data_sub} are generated with $S_j^{(m)}(0)$ randomly sampled from state variable space $I_S = [-2,2]$. The total number of data pairs is $1458$ in this case.

\begin{figure}[h]
\includegraphics[width = 1\textwidth]{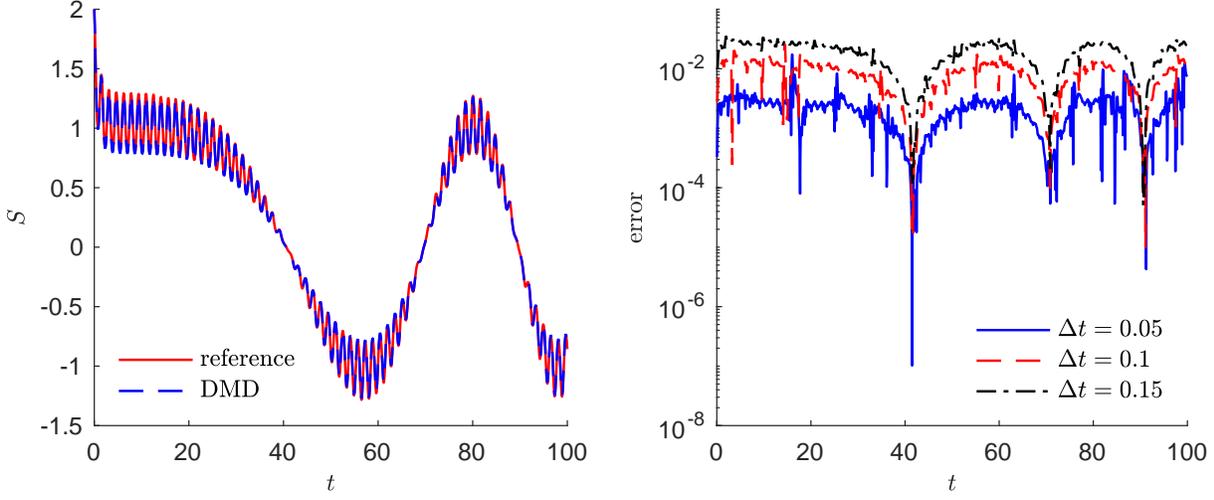}
\caption{(a). DMD prediction of~\eqref{eq:4-1} with inputs $\alpha(t) = \sin(4t)+1$ and $\beta(t) = \cos(t^2/1000)$; (b). The difference between the DMD surrogate and the reference solution using different $\Delta t$. }
\label{fig:4-1}
\end{figure}

Since $S$ is a scalar, the observable $\mathbf g$ can be designed to construct an augmented data matrix. Our choice of the observable is $\mathbf y (t_k;\mathbf p_k) = \mathbf g(S(t_k;\mathbf p_k)) = [S(t_k;\mathbf p_k), S(t_k;\mathbf p_k)^2]^\top$. Alternatively, one can construct the shift-stacked data matrices as introduced in~\cite{tu2013dynamic}. After the DMD surrogate model is trained, we use it to conduct system prediction. In Figure~\ref{fig:4-1}, the prediction result with a new initial condition $S_0 = 2$ and new external inputs $\alpha(t) = \sin(4t)+1$ and $\beta(t) = \cos(t^2/1000)$ is shown for time up to $T=100$. The reference solution is also shown for comparison. It can be seen that our surrogate model produces accurate prediction for this relatively long-term integration.

Figure~\ref{fig:4-1}(b) plots errors in predictions of~\eqref{eq:4-1} learned using different values of $\Delta t$, which exhibit a different behavior than DNN as reported in Figure 4.1 of~\cite{qin2021data}. Following the same analysis as conducted in Theorem 3.5 of~\cite{qin2021data}, the error of our method can be identified from two sources: 1. the approximation error of the parameterization, i.e., the difference between the true non-autonomous system~\eqref{eq:ODEs} and the modified system~\eqref{eq:NODE_modified}; 2. the error of the surrogate, i.e., the difference between the modified system~\eqref{eq:NODE_modified} and the surrogate model (\eqref{eq:ROM-iter} and~\eqref{eq:ROM-state} in our DMD method). The impacts of the two sources are well balanced in DNN so the overall error stays the same order for all choices of $\Delta t$. However, the second source of error is more dominant in DMD and thus the overall error shows sensitive dependence on the choice of $\Delta t$. It is well reported in the literature of DMD studies that the surrogate accuracy is highly dependent on the number of snapshots $N_\text{sanp}$, sampling frequency (i.e., $1/\Delta t$) and the underlying system itself. 

\subsection{Predator-prey model with control}
We now consider the following Lotka-Volterra predator-prey model with a time-dependent input $u(t)$:
\begin{equation}\label{eq:p-p}
\begin{aligned}
&\frac{dS_1}{dt} = S_1-S_1S_2+u(t),\\
&\frac{dS_2}{dt} = -S_2+S_1S_2.
\end{aligned}
\end{equation}

The local parameterization for the external input is conducted using quadratic polynomials, i.e., $\mathbf p\in \mathbb R^3$. We set $I_{\mathbf p} = [0,5]^3$ and the state variable space $I_{\mathbf S} = [0,5]^2$. Therefore, $N_\text{MC} = N_\text{traj} = 27$ parameter points on the Cartesian grid of $I_{\mathbf p}$ are selected to generate the training data with $N_\text{snap} = 9$ pairs. The total number of training data pairs is $243$. The observable is constructed in the following way:
$$\begin{aligned}
\mathbf y(t_k;\mathbf p_k) = \mathbf g(\mathbf S(t_k;\mathbf p_k)) =& [S_1(t_k;\mathbf p_k), S_2(t_k;\mathbf p_k), S_1(t_k;\mathbf p_k)^2, S_1(t_k;\mathbf p_k)S_2(t_k;\mathbf p_k),S_2(t_k;\mathbf p_k)^2,S_1(t_k;\mathbf p_k)^3, \\
&S_1(t_k;\mathbf p_k)^2S_2(t_k;\mathbf p_k), S_1(t_k;\mathbf p_k)S_2(t_k;\mathbf p_k)^2,S_2(t_k;\mathbf p_k)^3]^\top
\end{aligned}$$

\begin{figure}[h]
\includegraphics[width = 1\textwidth]{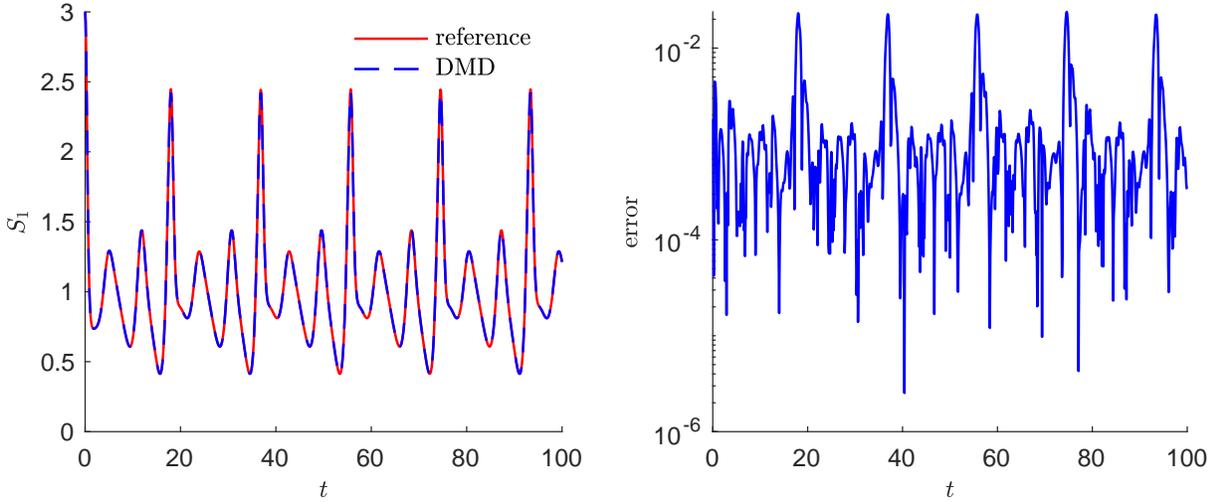}
\caption{DMD prediction of~\eqref{eq:p-p} with inputs $u(t) = \sin(t/3)+\cos(t)+2$ and $(S_1(0), S_2(0)) = (3,2)$. Left: $S_1$ solution profule; Right: Absolute error of $S_1$.}
\label{fig:4-2}
\end{figure}

The trained DMD surrogate is tested with external input $u(t) = \sin(t/3)+\cos(t)+2$ and initial condition $\mathbf S_0 = [3,2]$ for time up to $T = 100$. Figure~\ref{fig:4-2} shows the prediction profile and accuracy for $S_1$ compared with reference solution. The accuracy is satisfactory for this relatively long-time prediction. Similar results hold for $S_2$.

\subsection{Forced Oscillator}
We then consider a forced oscillator
\begin{equation}\label{eq:FO}
\begin{aligned}
&\frac{d S_1}{dt} = S_2,\\
&\frac{d S_2}{dt} = -\nu(t) S_1-S_2 +f(t),
\end{aligned}
\end{equation}
where the damping term $\nu(t)$ and the forcing $f(t)$ are time-dependent processes. Local parameterization for the inputs is conducted using quadratic polynomials and $N_\text{MC} = N_\text{traj} = 3^6$ parameter points on the Cartesian grid of $I_{\mathbf p} = [-1,1]^6$ are selected for training. On each parameter point, initial conditions are sampled randomly from state variable space $I_{\mathbf S} = [-3,3]^2$ and $N_\text{snap} = 3$ pairs of one-time step input-output data pairs are collected, making a total number of $2187$ pairs of training data.

\begin{figure}[h]
\includegraphics[width = 1\textwidth]{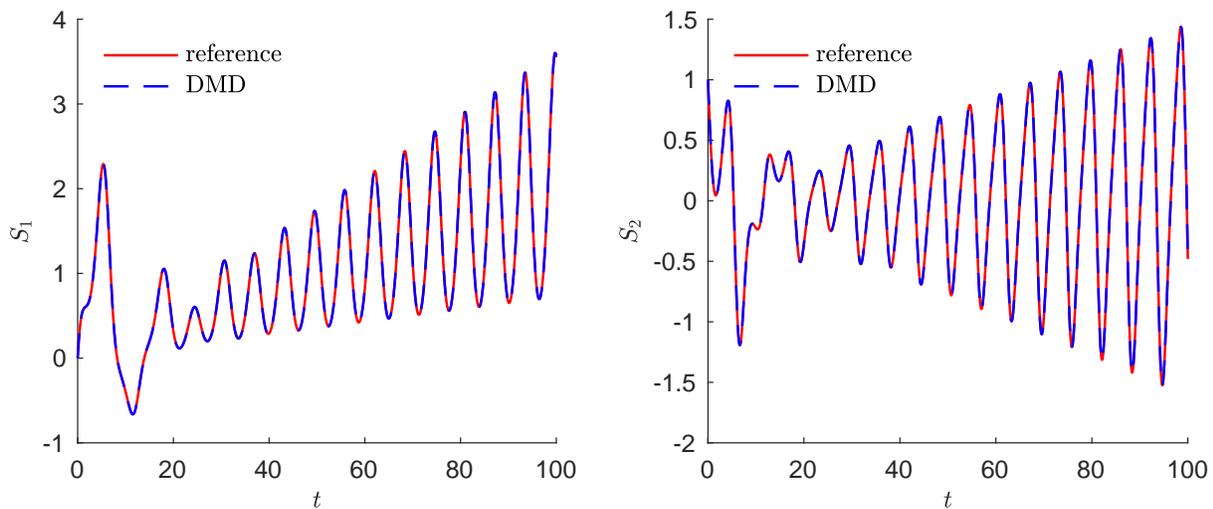}
\caption{DMD prediction of~\eqref{eq:FO} with inputs $\nu(t) = \cos(t)$ and $f(t) = t/200$.}
\label{fig:4-3}
\end{figure}

Predictions using the trained DMD model is shown in Figure~\ref{fig:4-3} for an arbitrarily chosen external inputs $\nu(t) = \cos(t)$ and $f(t) = t/200$. We observe very good agreement with the reference solution for relatively long-term simulation up to $T = 100$. Moreover, when $t>80$, the value of $x_1$ has been out of the training domain $I_{\mathbf S} = [-3,3]$, whereas the DMD model still generates accurate predictions. This shows some generalization ability of the DMD model.

%% number of data vs. DNN

\subsection{1D PDE: heat equation with source}
We now consider a PDE model, in particular, the following heat equation with a source term,
\begin{equation}\label{eq:1d_heat}
\begin{aligned}
&s_t = s_{xx}+q(t,x),\quad x\in [0,1], \\
&s(0,x) = s_0(x),\\
&s(t,0) =s(t,1) = 0,
\end{aligned}
\end{equation}
where $q(t,x)$ is the source term varying in both space and time. We set the source term to be
\begin{equation}
q(t,x) = \alpha(t)e^{-\frac{(x-\mu)^2}{\sigma^2}}
\end{equation}
where $\alpha(t)$ is its time-varying amplitude and parameters $\mu$ and $\sigma$ determine its spatial profile.

The learning of~\eqref{eq:1d_heat} is conducted in a discrete space. $N_S = 20$ equally distributed grid points are employed in the domain $[0,1]$. Let 
\begin{equation}\label{eq:4-4}
\mathbf S(t) = [s(t,x_1),\cdots, s(t,x_{N_{S}})]^\top.
\end{equation}
Local parameterization of the input $\alpha(t)$ is conducted using $4$th order polynomials. More specifically, $N_\text{MC} = N_\text{traj} = 3^7$ parameter points are selected from the Cartesian grid of the
local parameterization space $I_{\mathbf p} = I_{\alpha}\times I_{\mu}\times I_{\sigma} = [0,1]^5\times [0,3]\times [0.05,0.5]$ and $N_\text{snap} = 25$ training data pairs are generated by randomly sampling from state variable space $I_{\mathbf S} = [0,1]^{N_{S}}$ and doing one-time step simulation. There are $54675$ pairs of data in the training dataset. The we test the prediction ability of our surrogate model for a new source term (not in training data set), where $\alpha = t-\lfloor t \rfloor$ is a sawtooth discontinuous function, $\mu = 1$ and $\sigma = 0.5$.

\begin{figure}[h]
\includegraphics[width = 1\textwidth]{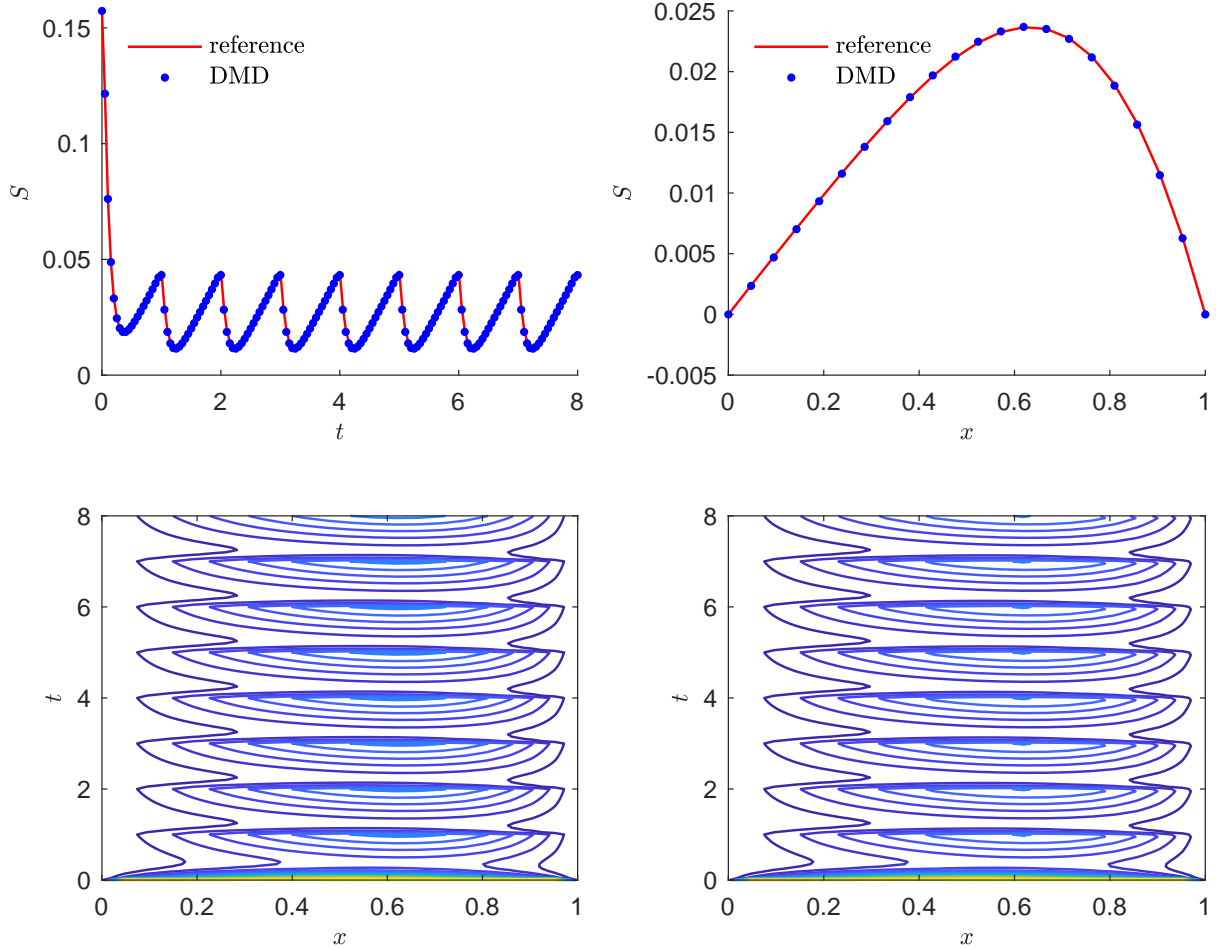}
\caption{DMD prediction of~\eqref{eq:1d_heat} with inputs $\alpha(t) = t-\lfloor t \rfloor$ , $\mu = 1$ and $\sigma = 0.5$. (a). Solution evolution at $x = 5$; (b). Solution profile at $t = 2$; (c). Reference solution contours over time; (d). DMD prediction contours over time.}
\label{fig:4-4}
\end{figure}

The prediction results are shown in Figure~\ref{fig:4-4} along with the reference solution solved from the underlying PDE. We observe satisfactory agreement between the surrogate model prediction to the reference solution. 

\section{Conclusion}
In this paper, we presented a numerical framework for learning unknown nonautonomous dynamical systems via DMD. To circumvent the numerical difficulties of computing the spectrum of the nonautonomous Koopman operator, the nonautonomous system is transformed into a family of modified  systems over a set of discrete time instances. The modified system, induced by a local parameterization of the external time-dependent inputs over each time instance, can be learned via our previous work dimension reduction and interpolation for parametric systems (DRIPS). The interpolation of the surrogate models in the parameter space allows one to conduct system predictions for other external time-dependent inputs by computing the parametric ROM of the new modified system iteratively over each time instance. The advantages of our proposed framework include: 1. Unlike previous work of approximating the spectrum of the time-depemdent Koopman operator,  our method works for general nonautonomous systems without any special requirements in special structures (e.g. periodic/quasi periodic). 2. Compared with other data-driven learning method like DNN, our method can achieve comparably satisfactory accuracy using much less training data. The efficiency and robustness of our method are demonstrated by various numerical examples of ODEs and PDEs. Future work along this line of research includes application of this framework to more complex large systems and rigorous analysis for the surrogate modeling errors.

\bibliography{DMD_nonautonomous}

\end{document}